\documentclass[11pt]{amsart}

\usepackage[T1]{fontenc}
\usepackage[utf8]{inputenc}
\usepackage[english]{babel} 
\usepackage{amsmath,amssymb,amsfonts,amsthm}
\usepackage{mathrsfs, esint, dsfont}
\usepackage{moreverb,rotating,graphics,float}
\usepackage{caption, subcaption}
\usepackage[colorlinks, linkcolor={blue}, citecolor={red}]{hyperref}
\usepackage{framed,fancybox}
\usepackage{enumerate}
\usepackage{slashed}
\usepackage{cancel}
\usepackage{color}
\usepackage{xcolor}
\usepackage[normalem]{ulem}
\usepackage{sidenotes}
\usepackage{setspace}

\usepackage{natbib}

\numberwithin{equation}{section} 
\newtheorem{theorem}{Theorem}[section]
\newtheorem{proposition}[theorem]{Proposition}
\newtheorem{remark}[theorem]{Remark}
\newtheorem{lemma}[theorem]{Lemma}
\newtheorem{corollary}[theorem]{Corollary}
\newtheorem{definition}[theorem]{Definition}

\newcommand{\ri}{\mathrm{i}}
\newcommand{\cS}{\mathscr{S}}

\newcommand{\cC}{\mathcal{C}}

\newcommand{\C}{\mathbb{C}}
\newcommand{\R}{\mathbb{R}}
\newcommand{\N}{\mathbb{N}}

\newcommand\bA{{\bold A}}

\newcommand{\bB}{\mathbf{B}}

\newcommand\bC{{\bold C}}

\newcommand\bh{{\bold h}}
\newcommand{\bH}{{\bold H}}

\newcommand\bm{{\bold m}}

\def\cB{{\mathcal B}}
\def\cC{{\mathcal C}}

\def\cH{{\mathcal H}}

\def\cS{{\mathcal S}}

\def\cU{{\mathcal U}}

\def\ra{{\mathrm{a}}}

\def\rb{{\mathrm{b}}}

\def\rd{{\mathrm{d}}}
\def\re{{\mathrm{e}}}
\def\ri{{\mathrm{i}}}
\def\rp{{\rm p}}

\def\rspan{{\mathrm{span}}}
\def\rsupp{{\mathrm{supp}}}

\newcommand{\wt}[1]{\widetilde{#1}}
\newcommand{\set}[1]{\left\{ #1\right\}}

\newcommand{\Ker}{{\rm Ker}}

\title[Spectrum of Iwatsuka type operators]{
On the Spectrum of the Iwatsuka Hamiltonian with Equal Magnetic Field Limits}

\author[A. Maichine]{Abdallah Maichine}
\address{Department of Mathematics, Faculty of Sciences, Mohammed V University in Rabat, 4 Av. Ibn Battouta, P.B. 1014 RP, Rabat, Morocco}
\email{a.maichine@um5r.ac.ma}

\date{\today}

\begin{document}
\begin{abstract}
We consider the two-dimensional magnetic Laplacian associated with a magnetic field that depends on a single spatial variable, known as the Iwatsuka Hamiltonian.
We focus on the equal-limit case, where the magnetic field is of the form $b(\cdot)=b_\infty+q(\cdot)$, with $b_\infty>0$ and $q\in C_0(\mathbb R)$. We prove that the Iwatsuka Hamiltonian admits an absolutely continuous spectrum whenever $q$ is compactly supported and has a constant sign near one edge of its support.
\end{abstract}

\maketitle
\tableofcontents

\section{Introduction}
In this paper, we investigate the spectral properties of the Iwatsuka-type magnetic differential operator, whose expression is
\begin{equation}\label{eq H_b introduction}
\cH_b :=-\partial_{xx}^2 +(-\ri\partial_y +\ra(x))^2,
\end{equation}
where, for every $x\in\R$, $a(x):=\int_0^x b(t) \, dt$ and $b\in\C(\R,\R)$ is a magnetic field intensity.

The Hamiltonian $\cH_b$ arises in the study of the two-dimensional motion of quantum particles under the action of a perpendicular magnetic field $\bB=(0,0,b(x))$.

For a constant magnetic field $b(x)=b_0>0$, $\cH_b$ reduces to the well-known two-dimensional Landau Hamiltonian, whose spectrum is given by $\sigma(\cH_b)=b_0 (2\N-1)$, the so-called Landau levels, each of which is an eigenvalue of infinite multiplicity. In this case, $\cH_b$ is invariant under magnetic translations, resulting in a uniform electron density, 
see for instance \cite{landau_quantum_1958} or \cite[Section 2]{G-C_L_M_25}.\\
Any variation of the magnetic field breaks this invariance along the direction in which the field varies. It is conjectured that such a variation spreads the Landau levels into purely absolutely continuous spectrum, see \cite[Section 6.5]{cycon1987schrodinger} or \cite[Section 7.2.3]{ExnKov15}, thereby allowing electron propagation.

To the best of our knowledge, this conjecture has not yet been solved in its full generality. Indeed,
interest into this question for Hamiltonian $\cH_b$ were initiated with A. Iwatsuka in his seminal work \cite{Iwatsuka85}. Since then, the spectral properties of this particular class of Hamiltonians have attracted considerable attention, see for example~\cite{ExnKov00,miranda2016AHP,mirandaPopov2018,ManPur97,Tus16}. 
Related transport properties of this model have been investigated in \cite{brietRaikovSoccorsi2008,hislop2015edge,DomGerRaik2011quantization}. 

A fundamental feature of the Iwatsuka Hamiltonian $\cH_b$ is that its spectral properties are strongly related to the asymptotic behaviour of the magnetic field at $\pm\infty$. In particular, when the limits
\[
b_\infty^-:=\lim_{x\to-\infty}b(x),
\qquad
b_\infty^+:=\lim_{x\to+\infty}b(x)
\]
exist and are different, the spectrum of $\cH_b$ can be shown to be absolutely continuous.

The situation is less straightforward when the two asymptotic limits coincide. In that case, the information provided by the asymptotic magnetic field alone does not yield such a conclusion. 

The aim of this work is to study the absolute continuity of the spectrum of $\cH_b$ in the equal-limit case. 
More precisely, we consider magnetic fields satisfying $$\lim_{-\infty} \rb(x)=\lim_{+\infty} \rb(x)=b_\infty>0.$$ 


To better understand this dichotomy, let us agree that, due to its translation invariance in the $y$-direction, $\cH_b$ can be seen unitarily equivalent to a direct integral of fibers as follows 
 \begin{equation}
     \cH_b \sim \int_\R^\oplus \bh_b(\xi)\, \rd \xi,
 \end{equation}
where, for every $\xi$, $\bh_b(\xi):=-\partial_{xx}^2+(\ra(x)-\xi)^2$.
In particular, 
$$\sigma(\cH_\rb)=\underset{\xi\in\R}{\cup}\sigma(\bh_\rb (\xi)).$$
Moreover, if $b$ satisfies $b_{\min} \le b(x) \le b_{\max}$, for some $b_{\max} \ge b_{\min} >0$, then $\bh_b(\xi)$ admits a compact resolvent and hence a discrete spectrum 
$$\sigma(\bh_b(\xi))=\{ \lambda_n(\xi): n\in\N \},\quad\forall \xi\in\R.$$
The maps $\lambda_n(\cdot)$, $n\in\N$, are known as \emph{dispersion functions}, and are responsible for the transport properties of the corresponding electronic system.

The spectral bands of $\bH_b$ are defined as $\bB_n:=\overline{\lambda_n(\R)}$, for every $n\in\N$, so that $\sigma(\cH_b)=\underset{n\in \N}{\cup}\, \bB_n$.

It was proved in~\cite[Proposition 2.2]{ManPur97} that 
$$\displaystyle\lim_{x\to\pm\infty} \rb(x)=\rb_\infty^\pm \implies\lim_{\xi\to\pm\infty}\lambda_n(\xi)=(2n-1)b_\infty^\pm, \quad\forall n\in\N. $$
This immediately yields
$$(2n-1)[\min(\rb_\infty^-,\rb_\infty^+),\max(\rb_\infty^-,\rb_\infty^+)]\subset \cB_n.$$ 
Consequently, \(\mathcal H_\rb\) has absolutely continuous spectrum in the unequal-limit case $\rb_\infty^-\neq \rb_\infty^+$. Furthermore, equality holds when the magnetic field \(b\) is monotone; see~\cite[Lemma 2.3 (ii)]{Iwatsuka85}. 
We have summarized the above properties in Proposition~\ref{prop bands in the different limits case} below.

The equal-limit case $\rb_\infty^- = \rb_\infty^+=b_\infty$, which is the focus of this paper, is considerably more delicate, since both asymptotic limits of each band function coincide at $(2n-1) b_\infty$.

In this paper, we address this critical situation of a field $b$ with equal limits at $\pm\infty$. More precisely, we make the following assumptions:
\begin{itemize}
    \item[\textbf{(A1)}] $b\in C(\R,\R)$ and there exist $b_{\max}\ge b_{\min}>0$ such that 
    $$\rb_{\min}\le \rb(x)\le \rb_{\max},\qquad \forall x\in\R.$$
    \item[\textbf{(A2)}] $b(x)=b_\infty+ q(x)$ with $q\in C_0(\R)$.
\end{itemize}
To the best of our knowledge, 
the absolute continuity of the spectrum has not been established for the full class of nonconstant magnetic fields satisfying \textbf{(A1)} and \textbf{(A2)}.

In the present paper, in addition to \textbf{(A1)} and \textbf{(A2)}, we prove the absolute continuity of the spectrum of \(\mathcal H_b\) under either of the following additional assumptions:
\begin{itemize}
    \item[\textbf{(A3)}] $q\in C_c(\R)$ and admits a constant sign near one of the edges of its support;
    \item[\textbf{(A4)}] $q$ admits exactly two nodal domains.
\end{itemize}
Let us emphasize that Assumption \textbf{(A3)} encompasses a broad class of compactly supported perturbations, including those considered in \cite[Theorem 2.6]{ExnKov00}. 
The latter result requires \(q\) to have a constant sign throughout its support, or to have a polynomial decay while maintaining a constant sign near the edges of its support.\\
Moreover, Assumption \textbf{(A3)} also covers perturbations $q$ with finite number of nodal domains, including those considered in \cite[Theorem 3.2]{ManPur97}. 
However, (A3) does not cover perturbations \(q\) whose nodal points\footnote{Here, nodal points refer to zeros of \(q\) across which the function changes sign, i.e., points separating two nodal domains.} accumulate at both edges of their support. Namely, if $q$ is of the form 
\begin{equation}\label{eq counter example q}
    q(x)= x(1-x)\sin\left(\frac{1}{x(1-x)}\right)1_{(0,1)}(x),\qquad\forall x\in \R.
\end{equation}

Our approach consists in proving that, for every $n\in\N$, the dispersion function $\lambda_n$ is nonconstant. To this end, we first establish a formula its derivative in Lemma~\ref{lm lambda'(xi) for q compact supp}. Using the Hellmann--Feynman formula, we then show that, for every $\mu\neq 0$, one has
\begin{equation}\label{eq intro lambda_n'}
\lambda_n^{'} (\xi) =\frac{2}{\mu}\int_{-\infty}^{+\infty} (b(x)-\mu)\left(\ra(x)-\xi\right)\psi_n^2(x,\xi),\quad\forall\xi\in\R,
\end{equation}
where $\psi_n(\cdot,\xi)$ denotes a normalized eigenfunction of $\bh_\rb(\xi)$ satisfying $$\bh_\rb(\xi)\psi_n(\cdot,\xi)=\lambda_n(\xi)\psi_n(\cdot,\xi). $$
If $q$ is compactly supported, we can choose $|\xi|$ sufficiently large so that $\rsupp(f)$ lies in the forbidden region for $\psi_n(\cdot,\xi)$.   
Taking $\mu=b_\infty$ in \eqref{eq intro lambda_n'}, and applying a Wentzel–Kramers–Brillouin (WKB) approximation to $\psi_n(x,\xi)$ in the forbidden region, we obtain the following asymptotic expansion for $\lambda_n^{'}$
\begin{equation}\label{eq introd asymptotic for lambda_n^'}
 \lambda_n^{'}(\xi)= \frac{2}{b_\infty}\,|\xi|\,\psi_n^2(0,\xi)\, \int_{\rsupp(q)} q(x)\re^{-\int_0^x \ra(t)\rd t}\re^{-2x|\xi|}(1+r_n(x,\xi))\,\rd x,   
\end{equation}
with $\|r_n(\cdot,\xi) \|_{\infty,\rsupp(q)}=O(|\xi|^{-1})$. Finally, using \eqref{eq introd asymptotic for lambda_n^'}, we show that \(\lambda_n\) is nonconstant whenever \textbf{(A3)} holds. 

Proposition~\ref{prop abs cont for 2 nodal domains} yields the result under Assumption \textbf{(A4)}, which is closely related to one of the assumptions considered in \cite{ManPur97}. For the sake of completeness, we provide a simpler proof in this setting.\\

The paper is organized as follows. In Section~\ref{sec elementary prop}, we introduce the general setting and establish some elementary spectral properties of the Iwatsuka Hamiltonian, including a description of its spectral bands. Section~\ref{sec main result} is devoted to the statement of our main results. Finally, in Section~\ref{sec proofs}, we establish some preliminary estimates and prove Theorem~\ref{thm absolute continuity of spectra} and Proposition~\ref{prop abs cont for 2 nodal domains}.
\section{General setting and elementary properties}\label{sec elementary prop}
Let $b\in\bC(\R,\R)$ be a magnetic field intensity satisfying 
\begin{enumerate}\label{assump A1}
\item[\textbf{(A1)}] There exist two positive constants $0<b_{\min} \le b_{\max} <\infty$ such that 
\[ b_{\min} \le b(x) \le b_{\max},\quad\forall x\in\R.
\]
\end{enumerate}
Define $\ra(x)=\int_0^x b(t)\, \rd t $, for all $x\in\R$, and consider the following differential expression
\begin{equation}
\cH_b := -\partial_{xx}^2 +(-\ri\partial_y +a(x))^2
\end{equation}
acting on smooth functions of $\R^2$.

The differential operator $\cH_b$ admits a non-negative self-adjoint realization in $L^2(\R^2,\C)$. This realization is associated to the (symmetric) sesquilinear form $\mathfrak{a}:\cC^\infty_c(\R^2,\C)\times \cC^\infty_c(\R^2,\C)\to\C$ given by
 
\[ \mathfrak{a}(f,g):=\int_{\R^2} \partial_x f\, \partial_x \overline{g}\; \rd x \rd y + \int_{\R^2} \left(-\ri \partial_y + \ra(x))f\right) \overline{\left(-\ri \partial_y + \ra(x)\right)g} \; \rd x \rd y,\quad 
\]
for all $f,g\in \cC^\infty_c(\R^2,\C).$\\
It is well-known that $\mathfrak{a}$ is closable and that the domain of its closed extension coincides with the magnetic Sobolev space $\bH_\bA^1(\R^2,\C)$, see for instance \cite[Chap. 7]{LL}, which is defined by
\[ \bH_\bA^1(\R^2,\C):=\set{ f\in L^2(\R^2,\C): (-\ri \nabla +\bA(x))f \in L^2(\R^2,\C)},
\] 
where $\bA(x,y)=(0,\ra(x))$ is a vector field on $\R^2$.
It is worth mentioning that $\bA$ acts as the (2d counter part of the) magnetic vector potential associated to the field $\bB(x,y,z)=(0,0, b(x))$.

\subsection{Spectral bands of $\cH_b$} Let us define $\cU: L^2(\R^2,\C)\to L^2(\R,L^2(\R,\C))$ as the unitary operator whose expression is given by
\begin{equation}
\left(\cU f\right)_\xi (x):=\left(\cU f (\xi) \right)(x)=\frac{1}{\sqrt{2\pi}}\int_\R \re^{+\ri \xi y} f(x,y) \, \rd y, \quad \forall f\in\cS(\R^2,\C),
\end{equation}
where $\cS(\R^2,\C)$ stands for the Schwartz space of smooth fast decaying functions of $\R^2$.

Therefore, we can write
\begin{equation}\label{eq integral decomp of H_b}
\cU \cH_b \cU^* = \int_\R^\oplus \bh_b(\xi) \, \rd \xi, 
\end{equation}
where, for all $\xi\in\R$, one has
\begin{equation}
\bh_b(\xi) := -\partial_{xx}^2 +( \ra(x)-\xi)^2.
\end{equation}
By assumption \textbf{(A1)}, one can see that $( \ra(x)-\xi)^2$ blows up at infinity for any $\xi\in\R$. So that, $\bh_b(\xi)$ has a compact resolvent, for all $\xi\in\R$. Therefore, $\bh_b(\xi)$ has a discrete spectrum consisting of eigenvalues only. Moreover, all eigenvalues are simple \cite{Iwatsuka85}.
 Let us denotes $\lambda_1(\xi) < \lambda_2(\xi)<...$ the eigenvalues of $\bh_b(\xi)$ in the increasing order.
 \begin{remark}\label{rmk construction of x_xi}
 It is also worth mentioning that, under assumption \textbf{(A1)}, $\ra$ is a continuous\footnote{$\ra$ is indeed $C^1$ function.} increasing function, since $\ra'(x)=b(x)\ge b_{\min}>0$. In addition, $\displaystyle\lim_{\pm \infty} \ra(x)=\pm \infty$. So that, $\ra$ is a continuous morphism of $\R$. Consequently, for every $\xi\in\R$, there exists a unique $x_\xi\in\R$ such that $\ra(x_\xi)=\xi$.
 \end{remark}
 \begin{proposition}
 Assume that \textbf{(A1)} is satisfied. One has
 \begin{itemize}
\item $\xi\mapsto\lambda_n(\xi)$ admits an analytic extension to an open set of $\C$ containing the real line $\R$.
\item 
 \begin{equation}
 \sigma(\cH_b) = \bigcup_{\xi\in\R}\sigma(\bh_b(\xi))=\bigcup_{n\in\N} \cB_n,
 \end{equation}
 where $\cB_n=\overline{\lambda_n(\R)}$ denotes the $n$-th band of $\cH_b$.

\item 
For all $\xi$ and $n$, there exists a normalized eigenfunction $\psi_n (\cdot,\xi)\in L^2(\R)$ such that
$$ \Ker\left(\bh_b(\xi)-\lambda_n(\xi)\right)=\rspan\set{\psi_n(\cdot,\xi)} . $$
\end{itemize} 
 \end{proposition}

A natural question that arises is whether $\cH_b$ has an absolutely continuous spectrum, which can be obtained if its bands $\cB_n$, $n\in\N$ are intervals of nonzero width. In the case of a constant magnetic field $b(x)=b_{\min}=b_{\max}$, it is well-known that $\cH_b$ has a point spectrum with infinite multiplicity for each eigenvalue.

\begin{proposition}\label{prop bands in the different limits case}
Assume \textbf{(A1)} is satisfied. Then,
\begin{enumerate}[i)]
\item For every $n\in\N$, one has
 \begin{equation}\label{eq B_n subset (2n-1)[b_min, b_max]}
 \cB_n\subseteq (2n-1)[b_{\min},b_{\max}].
 \end{equation}
    \item Assume that $\displaystyle\lim_{-\infty} \rb$ and  $\displaystyle\lim_{+\infty} \rb$ exist. Let us define
\[  \rb_\infty^-:=\min(\lim_{-\infty} \rb,\lim_{+\infty} \rb),\quad  \rb_\infty^+:=\max(\lim_{-\infty} \rb,\lim_{+\infty} \rb).\]
Then, 
\begin{equation}\label{eq B_n contains (2n-1)[b_inf^-,b_inf^+]}
    (2n-1)[b_\infty^-,b_\infty^+]\subseteq \bB_n.
\end{equation}
\end{enumerate}
\end{proposition}
As a consequence of the above proposition, it follows a description of the bands $\bB_n$, $n\in\N$, in the non equal-limit case, i.e., $\displaystyle\lim_{-\infty} b\neq\displaystyle\lim_{+\infty} \rb$. One has
\begin{corollary}
  Assume \textbf{(A1)} is satisfied. Assume that $\displaystyle\lim_{-\infty} \rb$ and  $\displaystyle\lim_{+\infty} \rb$ exist and are not equal. Then, $\cH_\rb$ admits an absolutely continuous spectrum with  finite number of gaps in the spectrum.  
\end{corollary}

\section{Main result}\label{sec main result}
In this section, in addition to Assumption \textbf{(A1)},  we will make the following assumption
\begin{enumerate}
\item[\textbf{(A2)}]\label{assump A2} $b(x)=b_\infty+q(x)$, with $q\in C_0(\R)$ and $b_\infty>0$.
\end{enumerate}
Here, $C_0(\R)$ refers to the space of continuous functions vanishing at infinity. Namely,
\[ C_0(\R)=\{ f\in C(\R,\R) : \lim_{|x|\to\infty} f(x)=0 \}. \]
 Our goal is to prove that the spectrum of $\cH_b$ is absolutely continuous for a large class of nonzero perturbation $q\in C_0(\R)$.

Before stating our main result, let us clarify some terminology that will be used in formulating our assumptions.
\begin{definition}\label{def nodal domains} 
Let $f:\R\to\R$ be continuous.
\begin{enumerate}
    \item We say that $f$ has two nodal domains if, and only if, there exists $p\in\R$ such that 
    \begin{itemize}
            \item $x\mapsto (x-p)f(x)$ does not change sign,
            \item there exist $x_-\in(-\infty,p)$ and $x_+\in(p,\infty)$ such that $f(x_-) f(x_+)\neq0$.
    \end{itemize}
\item Assume $f:\in C_c(\R)$ has a compact support. Let us denote 
\[ x_{\min} :=\inf\rsupp(f) \quad\text{and}\quad x_{\max} :=\sup\rsupp(f). \]
We say that $f$ admits  a constant sign near one edge of its support if, and only if, there exists $x_0\in(x_m,x_M)$ such that such that either
\[
f\, 1_{(x_{\min},x_0)}\geq 0,
\quad\text{or}\quad
f\, 1_{(x_{\min},x_0)}\leq 0,
\]
or the same property holds on $(x_0,x_{\max})$.
\end{enumerate}
\end{definition}
Now, we state the main result of this paper
\begin{theorem}\label{thm absolute continuity of spectra}
Assume \textbf{(A1)} and \textbf{(A2)} are satisfied. Assume, furthermore, that following is satisfied
\begin{enumerate}
\item[\textbf{(A3)}] $q\in C_c(\R)$ and $q$ has a constant sign near one edge of its support.
\end{enumerate}
Then, for every $n$, $|\bB_n|>0$. In particular, $\cH_b$ admits an absolutely continuous spectrum.
\end{theorem} 
The above Theorem~\ref{thm absolute continuity of spectra} covers previous results proving the same claim for a field that is constant outside a compact set, such as \cite[Theorem 2.6]{ExnKov00}, \cite[Assumption (B2)]{Iwatsuka85}, \cite[Theorem 3.2]{ManPur97}. 
Its proof is given in Section~\ref{sec proofs}.\\

Besides the compactly supported case, $\cH_b$ admits an absolutely continuous spectrum for $q$ having exactly two nodal domains. This statement can be deduced from \cite{ManPur97}. We state it here and provide a simple proof in Section~\ref{sec proofs} for the sake of completeness.
\begin{proposition}\label{prop abs cont for 2 nodal domains}
Assume \textbf{(A1)} and \textbf{(A2)} are satisfied. Assume that the following is satisfied
\begin{itemize} 
\item[\textbf{(A4)}] $q$ has two nodal domains.
\end{itemize}   
Then, $\cH_b$ admits an absolutely continuous spectrum.
\end{proposition}

\section{Proof of the main result}\label{sec proofs}
Before proceeding to the proof of Theorem~\ref{thm absolute continuity of spectra} and Proposition~\ref{prop abs cont for 2 nodal domains}, we are going to introduce useful elementary results for the eigenvalues $\psi_n(\cdot,\xi)$. 

Let us define
\begin{equation}\label{eq def of l_n(x,xi)}
\ell_n(x,\xi) := \left(\psi_n^{'}(x,\xi)\right)^2 - \left( (\ra(x)-\xi)^2-\lambda_n(\xi) \right)\psi_n^2(x,\xi),
\end{equation}  
and $\beta_n(x,\xi):=\frac{\partial\ell_n}{\partial x} (x,\xi)$. One has
\begin{equation}
\beta_n(x,\xi)=-2 b(x) (\ra(x)-\xi)\psi_n^2(x,\xi).
\end{equation}
From \cite[Lemma 3.4]{Iwatsuka85}, we infer the following properties of $\psi_n$, $\beta_n$ and $\ell_n$
\begin{proposition}\label{prop estimate on psi_n and so on}
Suppose that \textbf{(A1)} is satisfied. Let $\xi\in\R$ and $n\in\N$. Then,
\begin{itemize}
\item There exist positive constants $L_n$, $c_n$ and $\delta_n$ such that
\begin{equation}\label{eq estimate for psi_n}
 |x-x_\xi| > L_n \Longrightarrow  |\psi_n(x,\xi)|\le c_n \exp\left(-\frac{|x-x_\xi|}{\delta_n}\right),
\end{equation}
where $x_\xi:=\ra^{-1}(\xi)$, for every $\xi\in\R$.
\item $\ell_n(x,\xi)\to 0$ as $x\to\pm\infty$.
\end{itemize}
\end{proposition}
The crucial key tool in our proofs will be a formula for the derivative of $\lambda_n$, that we present in the Lemma below.
\begin{lemma}\label{lm lambda'(xi) for q compact supp}
Assume that \textbf{(A1)} and \textbf{(A2)} are satisfied, and let $n\in\N$. Then, for every $\mu\neq 0$, one has
\begin{equation}\label{eq lamnda_n^'(xi)=int q beta_n}
\lambda_n^{'} (\xi) =\frac{2}{\mu}\int_{-\infty}^{+\infty} (b(x)-\mu)\left(\ra(x)-\xi\right)\psi_n^2(x,\xi),
\end{equation}
for every $\xi\in\R$.
\end{lemma}
\begin{remark}
The constant $\mu$ in \eqref{eq lamnda_n^'(xi)=int q beta_n} is chosen independent of $\xi$. It is also possible to state the result making $\mu $ depending of $\xi$.     
\end{remark}
\noindent\textbf{Proof of Lemma~\ref{lm lambda'(xi) for q compact supp}.}
Let $n\in\N$ and $\xi\in\R$. One has
\[
  \lambda_n(\xi) = \langle \bh_b(\xi)\,\psi_n(\cdot\, ,\xi),\, \psi_n(\cdot\,,\xi) \rangle_{L^2(\R)}.
\]
Since $\|\psi_n(\cdot,\xi)\|^2_{L^2(\R)}=1$. Then, the \emph{Hellmann–Feynman formula} yields
\[\lambda_n^{'}(\xi):=
  \frac{\rd\lambda_n}{\rd\xi}(\xi)
  = \left\langle \frac{\rd\, \bh_b(\xi)}{\rd\xi}\,\psi_n(\cdot\,,\xi),\,\psi_n(\cdot\,,\xi) \right\rangle=-2\,\langle (\ra(x)-\xi)\,\psi_n,\psi_n \rangle.
\]

Hence,
\[
  \lambda_n'(\xi)  = -2\int_{\mathbb{R}} (\ra(x)-\xi)\,\psi_n^2(x,\xi)\,dx =\int_{\mathbb{R}} \frac{\beta_n(x,\xi)}{b(x)}\,dx.
\]
According to Proposition~\ref{prop estimate on psi_n and so on}, one can write
\begin{align*}
 \int_\R \beta_n(x,\xi)\, dx &=\lim_{s\to\infty}\int_{-s}^s \beta_n(x,\xi)\, dx \\
 &= \lim_{s\to\infty} \left(\ell_n(s,\xi) - \ell_n(-s,\xi) \right)\\
 &=0.    
\end{align*}
Therefore, a straightforward computation yields
\begin{align*}
  \lambda_n'(\xi)
  &= \int_{\R} \frac{\beta_n(x,\xi)}{b(x)}\,dx
     - \frac{1}{\mu}\int_\R \beta_n(x,\xi)\, dx \\
   &= -\frac{1}{\mu}\int_{\R} \beta_n(x,\xi)\,\frac{b(x)-\mu}{b(x)}\,dx\\
   &= \frac{2}{\mu}\int_{\R} (b(x)-\mu) \left( \ra(x)-\xi\right) \psi_n^2(x,\xi)\; \rd x.
\end{align*}

\subsection*{Proof of Theorem~\ref{thm absolute continuity of spectra}} Assume that \textbf{(A3)} is satisfied. Let $n\in\N$, it suffices to prove that $\lambda_n$ is not constant. 
To do so, we will rely on the formula~\eqref{eq lamnda_n^'(xi)=int q beta_n} for $\lambda_n^{'}$. Taking $\mu=b_\infty$ in \eqref{eq lamnda_n^'(xi)=int q beta_n}, one has  
\[ \lambda_n^{'} (\xi) =\frac{2}{b_\infty} \int_{\R} q(x) \left( \ra(x)-\xi\right) \psi_n^2(x,\xi)\; \rd x.
\]
Let $M>0$ is such that $\rsupp(q)\subset[-M,M]$. Then, 
\[\lambda_n^{'} (\xi) = \frac{2}{b_\infty} \int_{-M}^M q(x)\tilde{\beta}_n(x,\xi)\, dx,\]
where $\tilde{\beta}_n(x,\xi):=\left( \ra(x)-\xi\right) \psi_n^2(x,\xi)$, for all $x\in[-M,M]$ and $\xi\in\R$. 


We are going to prove that there exists $\xi\in\R$ such that $\lambda_n^{'} (\xi)\neq 0$. To this purpose, our approach will be based on an asymptotic expansion of $\tilde{\beta}_n(x,\xi)$, as $\xi\to\infty$, that is uniform in $x\in[-M,M]$.\\
Throughout this proof $c$ is referring to a positive constant that may change from line to line. It might depends on $n$, but not on $\xi$ and $x$.

We are going to split the proof in three steps, for the reader's convenience.\\
\underline{Step 1: A WKB formula for $\psi_n$.} Since $\psi_n$ solves the eigenvalue problem
\begin{equation}\label{eq eigenvalue problem psi_n}
- \partial_{xx}^2\psi_n(x,\xi) +(\ra(x)-\xi)^2 \psi_n(x,\xi)=\lambda_n(\xi)\psi_n(x,\xi)          
\end{equation}
and $\lambda_n(\xi)$ is the $n$-th eigenvalue (min-max principle) of the Schr\"odinger operator $\bh_b(\xi)$. Then, by  classical results on zeros of eigenfuctions to elliptic operators, see  \cite{Titchmarsh_eigenfunctions},
$\psi_n(\cdot,\xi)$ admits exactly $(n-1)$~zeros that are localized in the \emph{allowed region} $\Omega_n(\xi)$, where 
$$\Omega_n(\xi):=\{ x\in\R:(\ra(x)-\xi)^2\le \lambda_n(\xi)\}.$$
Taking into account \eqref{eq B_n subset (2n-1)[b_min, b_max]}, one gets 
$$\Omega_n(\xi)\subseteq \{ x\in\R:(\ra(x)-\xi)^2\le (2n-1)\rb_{\max}\}.$$
Furthermore, it is easy to see that 
$$|\ra(x)-\xi|\ge b_{\min}\,|x-x_\xi|,\quad \forall x\in\R. $$ Therefore, choosing $\xi$ such that 
\begin{equation}\label{eq dist(x_xi,[_M,M] > L_n}
\mathrm{dist}(x_\xi,[-M,M])>L_n:=\frac{\sqrt{(2n-1)b_{\max}}}{b_{\min}},    
\end{equation}
one obtains that $\Omega_n(\xi)\cap[-M,M]=\emptyset$. In particular, there exists a threshold $\mathfrak{t}_n>0$ such that $[-M,M]\subset \Omega_n(\xi)^c$, for every $\xi\in(-\infty,-\mathfrak{t}_n)$.\\
Now, let $\xi<-\mathfrak{t_n}$. Set $h=-\xi^{-1}$. Hence, \eqref{eq eigenvalue problem psi_n} becomes
\begin{equation*}
-h^2 \partial_{xx}^2\psi_n(x,-h^{-1}) +\left(1+h\,\ra(x)\right)^2 \psi_n(x,-h^{-1})=h^2\,\lambda_n(-h^{-1})\psi_n(x,-h^{-1}). 
\end{equation*}
Since, in this case, $[-M,M]\subset\Omega_n(\xi)^c$, we shall apply the semiclassical WKB approximation in the forbidden zone, see for instance \cite[Chapter 6]{olver1997asymptotics}, and gets
\begin{equation}\label{eq WKB approcimation}
    \psi_n(x,-h^{-1})= \psi_n(0,-h^{-1})\, \rp_h(x)^{-\frac{1}{2}}\, \exp\left(-h^{-1}\int_0^x \rp_h(t)\, dt\right)\left(1 +O_x (h)\right),
\end{equation}
for all $x\in[-M,M]$ and $h\to0$, where
\begin{equation}\label{eq expressopn of p_h(x)}
\rp_h(x):=\left( (1+h\ra(x))^2-h^2\lambda_n(-h^{-1}) \right)^{\frac{1}{2}},\qquad \forall x\in[-M,M],
\end{equation}
and, by boundedness of $\ra$ in $[-M,M]$, 
\[\sup_{x\in[-M,M]}|O_x(h)|\le c\,|h|.\]

\underline{Step 2. Asymptotic formula for $\lambda_n'$.}
From \eqref{eq expressopn of p_h(x)} and by boundedness of $\lambda_n$ and of $\lambda_n^{-1}$, we infer that 
\[
\rp_h(x)= 1+h\, \ra(x) + O(h^2),\qquad\forall x\in[-M,M].
\]
Since $\ra(\cdot)$ is bounded in $[-M,M]$, and $\lambda_n(\cdot)$ is bounded on $\R$, the above expansion is uniform in $x$. In addition, denoting $\tilde{\ra}(x)=\int_0^x \ra(t)\,\rd t$, for every $x\in[-M,M]$, one obtains
\[ h^{-1} \rp_h (x)=\frac{x}{h}+\tilde{\ra}(x) +O(h),\qquad \forall x\in[-M,M].. \]
One has also
\[ (\rp_h(x))^{-\frac{1}{2}}=1+O(h),\qquad \forall x\in[-M,M].  \]
It thus follow
\[
\psi_n(x,-h^{-1})= \psi_n(0,-h^{-1})\exp(-x\, h^{-1}-\tilde{\ra}(x)) (1+O(h)),\qquad \forall x\in[-M,M].
\]
Equivalently,
\[ \psi_n(x,\xi)= \psi_n(0,\xi)\, \re^{-\tilde{\ra}(x)}\,\re^{-x\,|\xi|} (1+O(\xi^{-1})),\qquad \forall x\in[-M,M],  \]
which yields
\begin{align*}
    \tilde{\beta}_n(x,\xi)&:=(\ra(x)-\xi)\psi_n^2(x,\xi)\\
    &=|\xi|\,\psi_n^2(0,\xi)\,\re^{-2x|\xi|-2 \tilde{\ra}(x)}(1+O(\xi^{-1}),\qquad
\end{align*}
for all $x\in[-M,M]$. Denoting $\tilde{q}(x)=q(x)\re^{-2\tilde{\ra}(x)}$, for every $x\in[-M,M]$, we then infer the following expansion for $\lambda_n^{'}$
\begin{equation}\label{eq lambda_n' asymptotic formula}
    \lambda_n^{'}(\xi)= \frac{2}{b_\infty}\,|\xi|\,\psi_n^2(0,\xi)\, \int_{-M}^M \tilde{q}(x)\re^{-2x|\xi|}(1+r_n(x,\xi))\,\rd x,
\end{equation}
with $\displaystyle\sup_{x\in[-M,M]}|r_n(x,\xi)|\le c\, |\xi|^{-1}$, for some positive constant $c>0$.

\underline{Step 3: Concluding arguments.} 
By Assumption \textbf{(A3)}, there exists $x_0\in(x_{\min},x_{\max})$ such that either $(x_{\min},x_0)$ or $(x_0,x_{\max})$ is a nodal domain for $q$, with $x_{\min}$ and $x_{\max}$ refers to the edges of $\rsupp(q)$ as in Definition~\ref{def nodal domains}. Assume that $(x_{\min},x_0)$ is a nodal domain for $q$. We expect that the reader will realise, while reading the arguments below, that the other case will be treated in a very similar way, so that our choice of nodal domain is without loss of generality.

Let us argue by contradiction and suppose that $\lambda_n^{'}\equiv0$. Let  $\xi\in(-\infty,-\mathfrak{t}_n)$. From \eqref{eq lambda_n' asymptotic formula} we infer that
\[
\int_{x_{\min}}^M \tilde{q}(x)\re^{-2(x-x_0)|\xi|}(1+r_n(x,\xi))\,\rd x=0.
\]
In particular,
\[
\int_{x_{\min}}^{x_0} \tilde{q}(x)\, \re^{-2|\xi|(x-x_0)}\, \rd x = \int_{x_{\min}}^M r_n(x,\xi) \,\tilde{q}(x)\, \re^{-2|\xi|(x-x_0)}\, \rd x - \int_{x_0}^M \tilde{q}(x)\, \re^{-2|\xi|(x-x_0)}\, \rd x 
\]
Taking into account the estimate for $r_n$, one obtains
\[
\left| \int_{x_{\min}}^{x_0} \tilde{q}(x)\, \re^{-2|\xi|(x-x_0)}\, \rd x \right| \le \frac{c}{|\xi|}\int_{x_{\min}}^{x_0} |\tilde{q}(x)|\, \re^{-2|\xi|(x-x_0)}\, \rd x +\frac{\|q\|_\infty}{|\xi|}.
\]
Therefore, for every $\xi\le\min(-c,-\mathfrak{t}_n)$,
\begin{align*}
0\le\left(1-\frac{c}{|\xi|}\right)\int_{x_{\min}}^{x_0} |\tilde{q}(x)|\, \rd x
&\le\left(1-\frac{c}{|\xi|}\right)\int_{x_{\min}}^{x_0} |\tilde{q}(x)|\, \re^{-2|\xi|(x-x_0)}\, \rd x\\
&\le \frac{c\,\|q\|_\infty}{|\xi|}.
\end{align*}
By letting $\xi\to-\infty$, one claims that $\int_{-M}^{x_0} |\tilde{q}(x)|\, \rd x=0$, which implies that $q=0$ on $(-M,x_0)$. This leads to a contradiction and completes the proof.\\
\subsection{Proof of Proposition~\ref{prop abs cont for 2 nodal domains}} Assume that \textbf{(A4)} is satisfied. Let us denote by $\alpha$ the nodal point of $q$. By Lemma~\ref{lm lambda'(xi) for q compact supp}, one has
\begin{align*}
\lambda_n^{'} (\xi) = &\frac{2}{b_\infty}\int_{-\infty}^{\alpha} q(x)\left(\ra(x)-\xi\right)\psi_n^2(x,\xi)\, \rd x\\
&+\frac{2}{b_\infty}\int_{\alpha}^{+\infty}q(x)\left(\ra(x)-\xi\right)\psi_n^2(x,\xi)\,  \rd x
\end{align*}
Now, taking $\xi=\ra(\alpha)$, one obtains that the integrand $ x\mapsto q(x) \left(\ra(x)-\xi\right) \psi_n^2(x,\xi)$ has the same sign along $\R$. Since both $\left(\ra(\cdot)-\xi\right) \psi_n^2(\cdot,\xi)$ and $q$ have a finite number of zeros,  $q(\cdot) \left(\ra(\cdot)-\xi\right) \psi_n^2(\cdot,\xi)\not\equiv0$. Therefore, $\lambda_n^{'}(\ra(\alpha))\neq0$, which ends the proof.

\subsection*{Acknowledgements} The author thanks his colleagues O. El-Fallah and M. Ziyat for insightful discussions.

\bibliographystyle{plain}

\bibliography{references}

@article{mirandaPopov2018,
  title={Spectrum of the Iwatsuka Hamiltonian at thresholds},
  author={Miranda, P. and Popoff, N.},
  journal={J. Math. Anal. Appl.},
  volume={460},
  number={2},
  pages={516--545},
  year={2018},
  publisher={Elsevier},
  doi={10.1016/j.jmaa.2017.11.056}
}

@article{miranda2016AHP,
  title={Eigenvalue asymptotics for a Schr{\"o}dinger operator with non-constant magnetic field along one direction},
  author={Miranda, P.},
  journal={Ann. H. Poincar{\'e}},
  volume={17},
  number={7},
  pages={1713--1736},
  year={2016},
  publisher={Springer},
  doi={10.1007/s00023-015-0445-6}
}

@book{olver1997asymptotics,
  title     = {Asymptotics and Special Functions},
  author    = {Olver, F. W. J.},
  year      = {1997},
  publisher = {A K Peters/CRC Press},
  address   = {Wellesley, Massachusetts},
  edition   = {2nd},
  isbn      = {978-1-56881-069-0}
}

@article{brietRaikovSoccorsi2008,
  title={Spectral properties of a magnetic quantum Hamiltonian on a strip},
  author={Briet, P. and Raikov, G. D and Soccorsi, E.},
  journal={Asym. Anal.},
  volume={58},
  number={3},
  pages={127--155},
  year={2008},
  publisher={IOS Press},
  doi={10.3233/ASY-2008-875}
}

@article{hislop2015edge,
  title={Edge states induced by {I}watsuka Hamiltonians with positive magnetic fields},
  author={Hislop, P. D and Soccorsi, E.},
  journal={J. Math. Anal. Appl.},
  volume={422},
  number={2},
  pages={1194--1224},
  year={2015},
  publisher={Elsevier},
  doi={10.1016/j.jmaa.2014.09.043}
}

@article{DomGerRaik2011quantization,
  title={Quantization of edge currents along magnetic barriers and magnetic guides},
  author={Dombrowski, N. and Germinet, F. and Raikov, G. D},
  journal={Ann. H. Poincar{\'e}},
  volume={12},
  number={6},
  pages={1169--1197},
  year={2011},
  publisher={Springer}
}

@book{cycon1987schrodinger,
  title={Schr{\"o}dinger Operators: With Application to Quantum Mechanics and Global Geometry},
  author={Cycon, Hans L. and Froese, Richard G. and Kirsch, Werner and Simon, Barry},
  series={Texts and Monographs in Physics},
  year={1987},
  publisher={Springer-Verlag},
  address={Berlin--Heidelberg--New York},
  isbn={978-3-540-16758-7}
}

@book{Titchmarsh_eigenfunctions,
    author =  {Titchmarsh, E. C.},
    title = {Eigenfunction Expansions Associated with Second-order Differential Equations},
    edition = {Second},
    publisher = {Oxford Univ. Press },
    year = {1962} 
}

@book{landau_quantum_1958,
  title = {Quantum Mechanics: Non-Relativistic Theory},
  author = {Landau, L. D. and Lifshitz, E. M.},
  translator = {Sykes, J. B. and Bell, J. S.},
  series = {Course of Theoretical Physics},
  volume = {3},
  publisher = {Pergamon Press},
  address = {London},
  year = {1958}
}

@article{Iwatsuka85,
    author = {Iwatsuka, A.},
    title = {Examples of absolutely continuous Schr\"odinger operators in magnetic fields},
    journal = {Publ. RIMS },
    year = {1985},
    volume = {21},
    pages = {385-401}
}

@article{ManPur97,
    author = {Mǎntoiu, M. and Purice, R.},
    title = {Some Propagation Properties of the {I}watsuka Model},
    journal = {Commun. Math. Phys.},
    volume ={188},
    pages = {691 – 708},
    year = {1997}
}

@article{ExnKov00,
    author = {Exner, P. and Kovarik, H.},
    title = {Magnetic strip waveguides},
    journal = { J. Phys. A: Math. Gen.},
    year ={2000},
    volume ={33},
    pages ={3297–3311}
}

@book{ExnKov15,
    author = {Exner, P. and Kovarik, H.},
    title = {{Q}uantum {W}aveguides},
    publisher = {Springer Cham},
    series = {{T}heoretical and {M}athematical {P}hysics},
    year = {2015}
}

@article{Tus16,
    author = {Tusek, M.},
    title = {On an extension of the {I}watsuka model},
    journal = { J. Phys. A: Math. Theor.},
    volume ={49},
    year = {2016},
    pages= {365205}
}

@article{G-C_L_M_25,
    author = {Garc\`ia-Cervera, C. and Lahbabi, S. and Maichine, A.},
    title = {Electronic structure models with 2D symmetries in the presence of magnetic fields},
    journal = {Ann. H. Poincaré},
    year = {2026},
    volume = {},
    url = {https://doi.org/10.1007/s00023-026-01715-x}
}

@book {LL,
    AUTHOR = {Lieb, E. H. and Loss, M.},
     TITLE = {Analysis},
    SERIES = {{G}raduate {S}tudies in {M}athematics},
    VOLUME = {14},
   EDITION = {Second},
 PUBLISHER = {American {M}athematical {S}ociety},
   ADDRESS = {Providence, RI},
      YEAR = {2001},
     PAGES = {xxii+346},
      ISBN = {0-8218-2783-9},
   MRCLASS = {00A05 (26-01 28-01 31-01 35J10 42-01)},
  MRNUMBER = {1817225 (2001i:00001)},
}

\end{document}